\documentstyle[12pt,twoside]{article}
\newtheorem{theorem}{Theorem}
\newtheorem{conjecture}{Conjecture}
\newtheorem{proposition}{Proposition}
\newtheorem{lemma}{Lemma}

\newtheorem{remark}{Remark}
\newtheorem{corollary}{Corollary}
\newtheorem{definition}{Definition}

\def\R{I\!\!R}
\def\C{I\!\!\!\!C}
\def\N{I\!\!N}
\def\Q{I\!\!\!Q}
\def\Z{I\!\!\!\!Z}
\def\vsni{\vskip 0.2cm}

\def\ui{[0,1]}

\def\B{{\cal B}}

\def\be{\begin{equation}}
\def\ee{\end{equation}}
\def\qed{\diamondsuit}
\def\O{\Omega}

\def\z{\zeta}

\def\s{\sigma}

\def\L{\cal{L}}
\def\P{\cal{P}}
\def\B{\cal{B}}

\def\M{\cal{M}}
\def\P{\cal{P}}

\def\di{\displaystyle}

\begin{document}

\title{On the spectrum of Farey and Gauss maps}
\date{\today}

\author{Stefano Isola \thanks{Dipartimento di Matematica e Fisica dell'Universit\`a
di Camerino and INFM, via Madonna delle Carceri, 62032 Camerino, Italy.
e-mail: $<$isola@campus.unicam.it$>$.}}
\date{}
\maketitle
\begin{abstract} 
\noindent
In this paper we introduce Hilbert spaces of holomorphic functions given by generalized Borel and Laplace transforms
which are left invariant
by the transfer operators of the Farey map and its induced transformation, the Gauss map, respectively.
By means of a suitable operator-valued power series we are able to 
study simultaneously the spectrum of both these operators
along with the analytic properties of 
associated dynamical zeta functions.
This construction establishes an explicit connection between
previously unrelated results of Mayer and Rugh (see \cite{Ma1} and \cite{Rug}).
\end{abstract}
\section{Introduction} 
The spectral analysis of transfer operators for smooth uniformly expanding maps of the unit interval $\ui$ is now
fairly well understood (see \cite{C}, \cite{Ba1}). 
The spectrum depends crucially on the function space considered
which is in general a Banach space. For 
Banach spaces of sufficiently
regular functions, e.g. the space ${\cal C}^k$ of $k$-times differentiable functions on $\ui$ with $k\geq 0$, the transfer operator is
{\sl quasi-compact}. This means that its spectrum is made out of a finite or at most countable set of isolated eigenvalues
with finite multiplicity (the discrete spectrum) and its complementary, the essential spectrum. 
The latter is a disk whose radius is a function both of $k$ and the expanding constant $\rho$ of the map 
(see e.g. \cite{CI}), in such a way that if we let $\rho \to 1$ from above (e.g. approaching an intermittency transition) 
the essential spectral radius 
tends to coincide with the spectral radius itself. 
In particular, in order to understand the nature of the spectrum lying under the `essential spectrum rug'
we have to consider increasingly smooth test functions as $\rho$ approaches $1$. 
This suggests, for instance, that for a type 1 intermittency model at the tangent bifurcation
point (see \cite{PM}) one should consider suitable spaces of analytic functions. 
In this paper we construct a Hilbert space ${\cal H}_0$ of analytic functions which is left
invariant by the transfer operator ${\cal P}$ of the Farey map (see below for definitions), a prototype of smooth
{\sl intermittent interval map},
having a neutral fixed point at the origin. As a result, the spectrum of ${\cal P}$ when acting
on ${\cal H}_0$ turns out to be the interval $[0,1]$ with embedded eigenvalues $0$ and $1$, plus a finite or countably infinite set of
eigenvalues of finite multiplicity. The latter is conjectured to be empty.
This would improve for this example a previous result obtained
by Rugh in a more general framework \cite{Rug}.
The above and related achievements are obtained by (a slightly modified version of) an inducing procedure which was introduced for the first time in the pioneering
study \cite{P1} (see also \cite{PS}, \cite{HI}, \cite{Is})
for a rather general class of intermittent interval maps. 
The main tool in this construction is an operator-valued function ${\cal Q}_z$ which enjoys simple algebraic
relations both with ${\cal P}$ and the transfer operator ${\cal Q}$ of the Gauss map, the latter being obtained by inducing the Farey 
map with respect to
the first passage time a subset of $\ui$ away from the neutral fixed point. The spectral properties 
of ${\cal Q}_z$ when acting on a Hilbert space ${\cal H}_1\subset {\cal H}_0$ are then suitably translated into
those of ${\cal Q}$ in ${\cal H}_1$ as well as ${\cal P}$ in ${\cal H}_0$.

\noindent
The paper is organized as follows:
Section \ref{P} is devoted to introduce the Farey-Gauss pair, 
briefly discussing some (mostly known) properties of these maps and of their invariant measures and ending with
a short account of their intimate connection
with number theory. Further material on these general facts can be found in \cite{Bi}, \cite{Ki}, \cite{F},
\cite{Ma2}. The main results are contained in the two subsequent
Sections. Section \ref{TO} deals with the spectral analysis of transfer  operators.
We first introduce the operator-valued function ${\cal Q}_z$ and establish simple algebraic identities 
(Proposition \ref{oprelation}). 
We then extend to ${\cal Q}_z$
some previous results of Mayer and Roepstorff (see \cite{MaR1}, \cite{MaR2}) for the Gauss transfer operator ${\cal Q}$ 
obtaining as a by-product 
an analytic continuation of ${\cal Q}_z$ outside the unit disk which is crucial to exploit 
the above identities for spectral analysis purposes (Proposition \ref{invspace}).
The main results on the spectrum of ${\cal P}$ (Theorems \ref{spectrump} and \ref{hank}) are then obtained by combining these identities with an explicit 
integral representation of ${\cal P}$ 
on the
Hilbert space ${\cal H}_0$ (Theorem \ref{resrep}).
In Section \ref{ZF} we apply the construction of the previous Section to study analytic properties of the 
dynamical zeta functions \cite{Ba2} for the Farey-Gauss pair. The role of ${\cal Q}_z$ is here
played by a two-variable zeta function $\z_2 (s,z)$ which simply relates to the Farey and Gauss zetas (Proposition \ref{twozeta}) and whose
analytic structure is directly connected to the spectrum of ${\cal Q}_z$ (Theorem \ref{twoanal}).
As a result, the zeta function of the Farey map
turns out to extend meromorphically to the cut plane $\C \setminus [1,\infty)$ (Corollary \ref{zetas}). 

\noindent
Finally, we point out that
some generalized version (involving a `temperature' parameter $\beta$) of these functions were previously studied in 
\cite{Ma1}, \cite{Ma2}, \cite{Ma3} for the Gauss map and in \cite{D} for the Farey map paired with an induced version 
conjugated to the Gauss map\footnote{
The inducing construction used in \cite{D}, the same as in \cite{P1}, is only slightly different from that used here
and several quantities, e.g. operators and zeta functions, dealt with there
are closely related to those discussed here. I thank one of the referees for having let me know about this work.}.
In the more general context of piecewise analytic map
with a neutral fixed point results yielding meromorphic continuation to the cut plane for zeta functions as well as regularized Fredholm determinants  
were obtained in \cite{Rug}.
 
\vsni
\noindent
{\bf Acknowledgment.} I thank the referees for the helpful criticism.
\vfill\eject
\section{Preliminaires}  \label{P}
We shall first consider the {\it Farey map} of the interval $[0,1]$ into itself defined as
\be
F(x)= \cases{  F_0(x),  &if $\; 0\leq x\leq 1/2$\, , \cr
               F_1(x), &if $\; 1/2 < x\leq 1$\, , \cr }
\ee
where
\be
F_0(x):={x\over 1-x}\quad\hbox{and}\quad F_1(x):=F_0(1-x)={1\over F_0(x)}={1-x\over x}\cdot
\ee
The inverse branches are
\begin{eqnarray}\label{invbran}
\Psi_0(x)&\equiv &F_0^{-1}(x)= {x\over 1+x}={1\over 2} -{1\over 2}\left({1-x\over 1+x}\right),\nonumber \\
\Psi_1(x) &\equiv& F_1^{-1}(x)={1\over 1+x}={1\over 2} +{1\over 2}\left({1-x\over 1+x}\right)\, .
\end{eqnarray}
For $x\ne 0$ 
the map $\Psi_0(x)$
is conjugated
to the right translation $x\to S(x)= x+1$, i.e.
\be\label{conju1}
\Psi_0 =J\circ S \circ J \quad\hbox{with}\quad 
J(x)=J^{-1}(x)=1/x.
\ee
This yields for the $n$-iterate
\be\label{iterates}
\Psi_0^n(x) = J\circ S^n \circ J (x) = {x\over 1+nx}\cdot 
\ee
Moreover $\Psi_1(x)$ satisfies
\be\label{conju2}
\Psi_1(x) = J\circ S (x).
\ee
\subsection{The induced map}
Let ${\cal A}=\{A_n\}_{n\geq 1}$ be the countable partition of $[0,1]$ given by
$A_n=[1/(n+1),1/n]$. Setting $A_0=[0,1]$ it is easy to check that $F(A_n)=A_{n-1}$ for
all $n\geq 1$. 
Let $X$ be the residual set of points in $[0,1]$ which are not
preimages of $1$ with respect to the map $F_0$, namely
$X = (0,1] \setminus \{1/n\}_{n\geq 1}$.
The {\it first passage time} $\tau : X \to \N$  
in the interval $A_1$ is defined as
\be\label{tau}
\tau (x)= 1+ \min \{n\geq 0 \;:\; F^n(x)\in A_1\;\}=\left[{1\over x}\right],
\ee
where $[a]$ is the integer part of $a$.
We see that $A_n$ is the closure of the set $\{x\in X\, :Ê\, \tau (x)=n\}$.
On the other hand, the {\it return time} function $r : A_1 \to \N \cup \{\infty\}$ 
in the interval $A_1$ is given by
\be\label{erretau}
r (x)=\min \{n\geq 1 \;: \; F^n(x)\in A_1\;\}=\tau \circ F(x).
\ee

\noindent
We now consider the map $G:X\to X$
obtained from $F$ by inducing w.r.t. the first passage time $\tau$, i.e.
\be
G(x) = F^{\tau (x)}(x),
\ee
which can be extended to all of $\ui$ setting
$G(0)=1$, $G(1)=0$,
$$
\lim_{x\uparrow 1/n}G(x)=0,\quad \lim_{x\downarrow 1/n}G(x)=1,\quad n> 1,
$$
and whenever $x\in (1/(n+1),1/n)$ we have, using (\ref{iterates}),
\be
G(x) \equiv G_n(x)=F^n(x)=F_1\circ F_0^{n-1}(x)={1\over x} -n= {\di 1\over \di x}-\tau(x).
\ee
In other words the induced map is the celebrated {\sl Gauss map}
\be 
G(x) =\cases{ \left\{{\di 1\over \di x}\right\},  &if $\; x\neq 0$ , \cr
             \;\;\, 0, &if $\; x=0$ , \cr }
\ee
where $\{a \}$ denotes the fractional part of $a$.
It has countably many inverse branches $\Phi_n$ given by
\be
\Phi_n(x) = G_n^{-1}(x) = {1\over x+n},\qquad n\geq 1.
\ee
\subsection{Invariant measures} 
It is an easy task to verify that the $\s$-finite 
absolutely continuous measure 
\be\label{sigmafinite}
\nu (dx) \equiv e(x)\, dx ={1\over \log 2}\cdot
{\di dx\over \di x}
\ee
is invariant for the dynamical system $(\ui , F)$.
Note that $\nu (A_n)= (\log 2)^{-1} \log{(1+{1\over n})}$ and $\nu ([0,1])=\infty$.
Let $B_n=\{x\in A_1\, : \, r(x)=n\}$.  
Using (\ref{erretau}) we have ${\overline {F_1(B_n)}}=A_n$.
We now show that
$\nu (A_n) = \sum_{k\geq n} \nu (B_k)$. Indeed, for $n=1$ we have
$\sum_{k\geq 1}\nu (B_k) = \nu (A_1)=1$. Moreover, since $\nu$ is
$F$-invariant,
$\nu (A_n)=\nu(F^{-1}(A_n))=\nu(A_{n+1})+\nu(B_{n+1})$, and the assertion
follows by induction.
Therefore the expected return time is infinite:
\be \label{infty}
\nu_{A_1} (r) = \int_{A_1}\, r(x)\, \nu (dx)= \sum_{n\geq 1} n \, \nu (B_n) =\sum_{n\geq 1} \nu (A_n) =
\nu (\ui) =\infty,
\ee
where $\nu_{A_1}$ is the conditional probability measure defined as $\nu_{A_1}(E) = \nu(E\cap A_1)/\nu (A_1)$.
It is known that in this situation there 
is the coexistence of two different statistics for
the dynamical system $(F,\ui)$: besides $\nu$, 
the ergodic means ${1\over n}\sum_{i=0}^{n-1}\delta_{F^i(x)}$
converge weakly to the Dirac delta at $0$ (see \cite{Me}, \cite{HY}).

\noindent
Let $\rho$ be the probability measure obtained by pushing forward $\nu$ with $F_1$, i.e. 
\be\label{mes}
\rho (E) =   ((F_{1})_*\,\nu) (E)=
(\nu \circ \Psi_1) (E) .
\ee
Reasoning as above one readily verifies that
the converse relation is
\be
\nu (E) = \sum_{n\geq 0} (\rho \circ \Psi_0^n)(E).
\ee
In particular we have $\nu (A_n) = \sum_{l\geq n} \rho ( A_l)$ and 
$\rho (A_n) = \rho (F_1(B_n))= \nu (B_n)$, where $B_n$ is as above.
We then have
\be
\rho (E) = (\nu \circ \Psi_1)(E)
=\sum_{n\geq 0} (\rho \circ \Psi_0^n\circ \Psi_1)(E)= \rho (G^{-1}E),
\ee
which says that $\rho$ is $G$-invariant. Moreover $\rho$ is ergodic with respect to $G$ (see e.g. \cite{Bi}).
Setting $h(x) = \rho (dx) /dx$ we get
\be\label{densities}
h = |\Psi_1'|\cdot e\circ \Psi_1,\qquad
e = \sum_{k=0}^{\infty}(\Psi_0^{k})^{\prime} \cdot h\circ 
\psi_0^{k},
\ee
which gives the well known result
\be\label{gaussdensity}
h(x) =  {1\over \log 2}\cdot
{ dx\over (1+ x)}\, \cdot
\ee
The primitive $H(x)$ of $h(x)$, with $H(0)=0$, is
$H(x) =  { \log (1+x) / \log 2 }$.
Setting $q_n := H({1\over n+1})=(\log 2)^{-1} \log{(1+{1\over n+1})}$, 
we have
$\nu(A_n)=q_n$ and $\rho (A_n)=q_{n-1}-q_n$.
We see that $q_n$ is a (strict) Kaluza sequence, i.e. for all $n\geq 1$
\be
1 = q_0 > q_1 > \cdots > q_n >0\qquad\hbox{and}\qquad  q_n^2 < q_{n-1}\, q_{n+1}.
\ee

\noindent
Finally, by (\ref{tau}), (\ref{erretau}), (\ref{infty}) and (\ref{mes}) we have
\be\label{infty2}
\rho (\tau) = ((F_{1})_*\,\nu)(\tau) = \nu (\tau \circ F_1) = \nu (r) =\infty.
\ee
On the other hand we have the following,
\begin{lemma}\label{K}
 The function $\log \tau$ is in $L_1(\rho)$ 
and satisfies
\be \label{finite}
\lim_{n\to \infty} {1\over n}\sum_{j=0}^{n-1}\log \tau(G^j(x))= \rho(\log \tau) = K,\qquad \rho-{\rm a.e.},
\ee
where the positive constant $K$ is defined by
\be
e^K = \prod_{k=1}^\infty \left(1+{1\over k(k+2)}\right)^{\log k \over \log 2}.
\ee
\end{lemma}
{\sl Proof.} We have
\begin{eqnarray}
\rho(\log \tau) &=& \sum_{k=1}^\infty   \rho (A_k)\cdot \log k = \sum_{k=1}^\infty (q_{k-1}-q_k)\cdot\log k \nonumber
\\
&=& \sum_{k=1}^\infty {\log k \over \log 2}\cdot \log\left( \left(1+{1\over k}\right)\left(1+{1\over k+1}\right)^{-1}\right)\nonumber
\\
&=& \sum_{k=1}^\infty {\log k \over \log 2}\cdot \log\left( 1+{1\over k(k+2)}\right) = K <\infty.\nonumber
\end{eqnarray}
This computation shows both that $\log \tau\in L_1(\rho)$ and the last equality in (\ref{finite}).
The first equality in (\ref{finite}) now follows from the ergodic theorem \cite{Bi}. $\qed$
\vsni
\noindent
The constant $K$ which appears above is known in number theory as {\sl Khinchin's constant}. 
This is not a coincidence,
as we now briefly explain.
\subsection{Connection with number theory}  The {\rm Farey sum} over two rationals ${a\over b}$ and ${a'\over b'}$ is
the mediant operation given by \cite{HR}
\be
{a''\over b''} = {a+ a'\over b+b'}\cdot
\ee 
It is easy to see that ${a''\over b''}$ falls in the interval $({a\over b},{a'\over b'})$.
Now, having fixed $n\geq 0$, let ${\cal F}_n$ be the ascending sequence
of irreducible fractions between $0$ and $1$ obtained inductively in the following way.
Set first ${\cal F}_0=({0\over 1}, {1\over 1})$. Then ${\cal F}_n$ is obtained from ${\cal F}_{n-1}$
by inserting among each pair of consecutive rationals ${a\over b}$ and ${a'\over b'}$ in ${\cal F}_{n-1}$
their mediant ${a''\over b''}$ as above. Thus
${\cal F}_1=({0\over 1},  {1\over 2}, {1\over 1})$,
${\cal F}_2=({0\over 1},  {1\over 3}, {1\over 2},  {2\over 3}, {1\over 1})$, 
${\cal F}_3=({0\over 1}, {1\over 4}, {1\over 3}, {2\over 5},
 {1\over 2}, {3\over 5}, {2\over 3}, {3\over 4}, {1\over 1})$ and so on. The elements
of ${\cal F}_{n}$ are called {\rm Farey fractions}. The
name of the map $F$ can be related to the easily verified observation that the set of pre-images 
$\cup_{k=0}^{n+1} F^{-k}\{0\}$ coincides with ${\cal F}_n$ for all $n\geq 0$.
In particular, this implies that 
$\cup_{k=0}^\infty F^{-k}\{0\}=\Q \cap [0,1]$ (notice that the same is true for the induced map:
$\cup_{k=0}^\infty G^{-k}\{0\}=\Q \cap [0,1]$).

\noindent
On the other hand, we recall that
every real number $0<x<1$ has a continued
fraction expansion of the form \cite{Ki}
\be
x= {1\over \displaystyle k_1 + {1\over \displaystyle k_2 + {1\over\displaystyle k_3 +\cdots }}}=
[k_1,k_2,k_3,\dots]\, ,
\ee
with $k_i\in \N$. 
By applying Euclid's algorithm one sees that the above expansion terminates if and only if $x$ is a
rational number. There is an intimate connection between the partial quotients $k_1,k_2,\cdots$ and
the Gauss map $G$. Indeed, given $x$ as above we can write
\begin{eqnarray}
x={1\over {\di 1\over \di x} }&=&{1\over \left[{\di 1\over \di x}\right] + \left\{ {\di 1\over \di x}\right\} }
={1\over k_1 + G(x) }={1\over \displaystyle k_1 + {1\over {\di 1\over \di G(x)}} }\nonumber \\
&=&{1\over \displaystyle k_1 + {1\over \displaystyle  \left[{\di 1\over \di G(x)}\right] +
\left\{ {\di 1\over \di G(x)}\right\}   } }={1\over \displaystyle k_1 + {1\over \displaystyle  k_2 +
G^2(x)   } }=\cdots 
\end{eqnarray}
Therefore, $k_1=[1/x]$, $k_2=[1/G(x)]$, $k_3=[1/G^2(x)]$ and so on.
Alternatively,
\be\label{Gaction}
\hbox{if}\quad x=[k_1,k_2,k_3,\dots]\quad\hbox{then}\quad G(x) =[k_2,k_3,\dots]\, .
\ee
Farey fractions have close relationships with continued fractions.
Let us say that a Farey
fraction has order $n$ if it belongs to ${\cal F}_n\setminus {\cal F}_{n-1}$. Given $n\geq 1$ there are exactly
$2^{n-1}$ Farey fractions of order $n$ (they form the set $F^{-(n+1)}\{0\}$) and
it is possible to show (see below eq. (\ref{Faction})) that the integers $k_i$ 
in their (finite)
continued fraction expansion sum up to $n+1$. Furthermore, it is easy to realize that all
Farey fractions which fall in the interval $(1/(n+1),1/n)$ have order greater than or equal to $n+1$,
whereas their
continued fraction expansion starts with $k_1=n$. Thus, the map
$F$ acts on Farey fractions by reducing their order of one unit. We can write an explicit
expression for the action of $F$
on continued fraction expansions. Indeed, if $1/2<x\leq 1$ then $k_1=1$ and $F(x)={1\over x}-k_1= G(x)$.
If instead $0<x\leq 1/2$ then $k_1>1$ and $F(x)=1/({1\over x}-1)$. Therefore,
\be\label{Faction}
\hbox{if}\quad x=[k_1,k_2,k_3,\dots]\quad\hbox{then}\quad F(x) =[k_1-1,k_2,k_3,\dots]\, ,
\ee
with $[0,k_2,k_3,\dots]\equiv [k_2,k_3,\dots]$
(compare to (\ref{Gaction})).
Now, it is well known that for almost all $x\in (0,1)$ the arithmetic mean of the partial quotients is infinite
(see, e.g., \cite{Ki}),
i.e.
\be 
\lim_{n\to \infty} {k_1 +\cdots +k_n \over n} = \infty,\qquad\hbox{(a.e.)}
\ee
From the above discussion and (\ref{tau}) we get $k_l=[1/G^{l-1}(x)]=\tau (G^{l-1}(x))$, which for $l>1$ is
the time between the $(l-1)$-st and the $l$-th passage in $A_1$ of the orbit of $x$ with $F$. 
Therefore,
the total number $S_n$ of iterates of $F$ needed to observe $n$ passages 
in $A_1$, that is the function
\be
S_n(x)=\tau (x)+\tau(G(x))\cdots +\tau(G^{n-1}(x)),
\ee
satisfies
\be \label{et}
\lim_{n\to \infty} {S_n(x) \over n} = \infty\qquad\hbox{(a.e.)}
\ee
Since $\rho$ is absolutely continuous w.r.t. the Lebesgue measure on $[0,1]$, the properties expressed by (\ref{infty2}) and (\ref{et})
can be regarded as an instance of validity of the ergodic theorem for the non-integrable
function $\tau$.
One can actually say more. As a consequence of (\cite{Ki}, Theorem 30) we have that for almost all $x\in (0,1)$  
the inequality
\be
S_n(x) \geq n\, \log n
\ee
is satisfied for an infinite number of values of $n$. On the other hand, Lemma \ref{K} can now be rephrased
by saying that the geometric mean 
of the partial quotients has a certain finite value (a.e.).
This, in turn, is a corollary
of a theorem of Khinchin (\cite{Ki}, Theorem 35), which says that for any function $f(k)$ defined on the positive integers
and satisfying $f(k)={\cal O}(k^p)$ with $0\leq p< 1/2$ we have, for almost all $x\in (0,1)$,
\be
\left| {1\over n}\sum_{j=i}^{n}f(k_j) - \sum_{k=1}^\infty {f( k) \over \log 2}\cdot \log\left( 1+{1\over k(k+2)}\right)\right| \leq \epsilon (n)
\ee
where the error function $\epsilon (n)$ is any positive function decreasing to zero as $n\to \infty$ so that
$\sum  n^{-2}\cdot \epsilon^{-2}(n) < \infty$. Lemma \ref{K} then  corresponds to the choice $f(k)=\log k$.

\section{Transfer operators}\label{TO}
We start by establishing some formal algebraic relations between
the transfer operators $\P$ and $\M$ associated to the maps $F$ and $G$, 
respectively (see \cite{Ba1}).
They describe the action of the differentiable dynamical systems $F$ and $G$ on the density 
$f$ of a measure absolutely continuous measure wrt Lebesgue by
\begin{eqnarray}\label{opP}
{\P} f (x ) &=& ({\P}_0 + {\P}_1 ) f (x)
=: |\Psi_0'(x)|\cdot  f(\Psi_0(x))+|\Psi_1'(x)|\cdot f(\Psi_1(x) )\nonumber \\
&=& \left({1\over x+1}\right)^2 \left[ f\left( {x\over x+1}\right) + f\left( {1\over x+1}\right) \right],
\end{eqnarray}
and
\be\label{emmezeta}
{\cal Q} \, f (x) =
\sum_{n=1}^{\infty} Q_n \, u(x)=:\sum_{n=1}^{\infty} |\Phi'_n(x)| \cdot f (\Phi_n (x))
=\sum_{n=1}^{\infty}\left({1\over x+n}\right)^2 \, f\left( {1\over x+n}\right).
\ee
We first notice that 
\be\label{rel1}
Q_n f(x) = {\P}^n (f\cdot \chi_n) (x)={\P}_1{\P}_0^{n-1}f(x),
\ee
where $\chi_n $ is the indicator function of $A_n$. Let ${\cal S}f(x) := f\circ S (x)=f(x+1)$ be the
shift operator. Note by (\ref{conju1}) and (\ref{conju2}) we have
\be\label{opconju}
 {\P}_1\,{\P}_0 f(x) = {\cal S} \, {\P}_1\, f(x),
\ee
and therefore (\ref{rel1}) yields
\be\label{doubleid}
{Q}_n f(x) = {\P}_1{\P}_0^{n-1}f(x)={\cal S}^{n-1}\,  {\P}_1 \, f(x).
\ee
More generally, for $z\in \C$, we shall consider a formal operator-valued power series
${\cal Q}_z$ defined by
\be\label{opevalpose}
{\cal Q}_z f(x) =\sum_{n=1}^{\infty} z^{\tau(\Phi_n (x))} \cdot|\Phi'_n(x)|\cdot f (\Phi_n (x)) 
=z\, {\P}_1(1-z{\P}_0)^{-1} f(x)
\ee
so that ${\cal Q}_1\equiv {\cal Q}$. 
The following operator relations are in force and are
independent of the function space the operators are acting on.
\begin{proposition}\label{oprelation} Let $z\in \C$ be such that (\ref{opevalpose}) 
is absolutely convergent. Then we have
\be\label{first}
(1-{\cal Q}_z)(1-z{\P}_0)=1-z{\P}
\ee
and
\be\label{second}
(1-z\, {\cal S})(1-{\cal Q}_z)=1-z\, {\tilde {\P}}\, .
\ee
where ${\tilde {\P}}={\cal S}+{\P}_1$.
\end{proposition}
\begin{remark} As already remarked in the Introduction,
an inducing procedure closely related to that
used here and leading to the study of the operator-valued function ${\cal M}_z =
(1-z{\P}_0)^{-1} z{\P}_1$ has been introduced in Prellberg's thesis \cite{P1} (see also \cite{PS} and \cite{D}), 
where algebraic identities
closely related to those stated above have been used to 
achieve a deep understanding of the thermodynamic formalism for intermittent interval maps.
Notice that ${\P}_1\, {\cal M}_z = {\cal Q}_z\, {\P}_1$.
Similar constructions have been used in \cite{HI} and \cite{Is}.
\end{remark}
\noindent
{\it Proof of Proposition \ref{oprelation}.} 
Using the first identity in (\ref{doubleid}) we get
\begin{eqnarray}
&&(1-{\cal Q}_z)(1-z{\P}_0) = (1- \sum_{n=1}^{\infty}z^n{\P}_1{\P}_0^{n-1})(1-z{\P}_0)= \nonumber \\
&&1-z{\P}_0-\sum_{n=1}^{\infty}z^n{\P}_1{\P}_0^{n-1}+\sum_{n=1}^{\infty}z^{n+1}{\P}_1{\P}_0^{n}
=1-z{\P}_0-z{\P}_1=1-z{\P}. \nonumber
\end{eqnarray}
In a similar way, using the second identity in (\ref{doubleid}) one shows (\ref{second}). $\qed$

\begin{corollary}\label{oprel2}
 Let $z\ne 0$ be such that (\ref{opevalpose}) is absolutely convergent and assume that the kernel of $1-z{\P}_0$ is empty. Then
$1$ is an eigenvalue of ${\cal Q}_z$ 
if and only if $z^{-1}$ is an eigenvalue both of
${\P}$ and ${\tilde {\P}}$, and they have the same geometric 
multiplicity. 
Furthermore, the
corresponding eigenfunctions $e_z$ of $\P$ and $h_z$ of ${\tilde {\P}}$ and 
${\cal Q}_z$ are related by $h_z = (1-z{\P}_0) e_z$ or else  
$e_z = \sum_{k=0}^{\infty} z^k{\P}_0^k h_z$.
\end{corollary}
\vskip 0.2cm
{\it Proof.} Assume that ${\cal Q}_zh_z = h_z$. From (\ref{first}) it then follows that 
$(1-z{\P})\sum_{k=0}^{\infty} z^k{\P}_0^k h_z = 0$. Conversely, 
assume that $z{\P} e_z = e_z$, then we have
$(1-{\cal Q}_z)(1-z{\P}_0)e_z=0$. In the same way, from (\ref{second}) it follows that ${\cal Q}_zh_z = h_z$ if and only if
${\tilde {\P}}h_z= z^{-1}h_z$.
$\qed$
\begin{remark}
As it will be clear in the sequel the condition on the emptyness of the kernel of $1-z{\P}_0$ is plainly satisfied
in the function space ${\cal H}_0$ considered below (cf (\ref{pizero})).
\end{remark} 
\begin{remark}
Setting $z=1$ in Proposition \ref{opevalpose} we recover (\ref{densities}) with
$e\equiv e_1$ and $h\equiv h_1$. In particular we see that 
the Gauss probability density (\ref{gaussdensity}) is a fixed point both of ${\cal Q}$ and 
${\tilde {\P}}$.
\end{remark} 
Having fixed an open connected domain $\O\subset \C$ let ${\cal H}(\O)$ be the Fr\'echet space
of functions which are holomorphic in $\O$ with the topology generated 
by the family of sup norms on compact subsets of $\O$. Moreover, we let 
$A_{\infty}(\O)\subset {\cal H}(\O)$ denote the Banach space given by
the subset of functions in ${\cal H}(\O)$ having continuous extension to ${\overline \O}$, endowed 
with the norm
\be
\Vert f\Vert = \sup_{w\in {\overline \O}} | f(w)|,
\ee
(where $w=x+iy$). We let first
${\cal Q}_z$ act on the Banach space $A_\infty(D)$ with
$D=\{w\in \C\, :Ê\, |w-1|<1\}$.
It is easy to verify that $\Phi_n ({\overline D}) \subset D$
for all $n\in \N$. Standard arguments (see \cite{Ma2}) then imply that whenever the power series
in (\ref{opevalpose})
is uniformly convergent ${\cal Q}_z $ defines a nuclear operator of order zero on $A_\infty(D)$.
\begin{lemma} The power series of 
${\cal Q}_{z}: A_\infty(D)\to A_\infty(D)$ has radius of convergence
bounded from below by $1$ and, moreover, 
it converges absolutely at every point of the unit circle.
\end{lemma}
{\it Proof.} 
The radius of convergence of ${\cal Q}_{z}$ is  
$\lim_{n\to \infty}\Vert Q_{n} \Vert^{-1/n}$ (here $\Vert \; \Vert$ denotes
the operator norm as well). We have
$\sup_{w\in {\overline D}}|Q_{n} f(w)| \leq C \, n^{-2}\,  \Vert f \Vert$
and therefore $\Vert Q_{n} \Vert\leq C \,  n^{-2}$. $\qed$
\vsni
\noindent 
We now introduce a subspace of $A_\infty(D)$ on which the action of ${\cal Q}_{z}$ will turn out to be
particularly expressive.
This is achieved via a generalized Laplace transform. 
\begin{definition}
We let
${\cal H}_1$ denote the Hilbert space of all complex-valued functions $f$ 
 which have a representation as generalized Laplace transform
\be\label{represent}
f(w) = ({\L}\, [\varphi])(w):=\int_0^\infty \, e^{-tw}\, \varphi (t)\, dm(t) 
\ee
where $\varphi \in L_2(m)$ and
$dm$ is the measure on $\R^+$ given by
\be\label{measure}
dm(t) = {t\, \over e^t-1}\, dt\, .
\ee
As a Hilbert space ${\cal H}_1$ is endowed with the inner product
\be
(f_1,f_2) = \int_0^\infty {\overline {\varphi_1(t)}}\, \varphi_2  (t)\, dm(t)\quad\hbox{if}\quad f_i={\L}\, [\varphi_i].
\ee
\end{definition}
\begin{remark}\label{h}
Putting $z=1$ we see that the $G$-invariant density $h$ can be
represented as $h=(\log 2)^{-1}\,{\L}[ \, (1-e^{-t})/t\,]$. 
\end{remark} 
The following Proposition generalizes
corresponding results obtained by Mayer and Roepstorff (see \cite{MaR1}, \cite{MaR2}) for the operator $\cal Q$.
\begin{proposition}\label{invspace}
 For each $z\ne 0$ with $|z|\leq 1$, the space ${\cal H}_1$
is invariant under ${\cal Q}_{z}$. More precisely we have
\be\label{rep}
{\cal Q}_{z}\, {\L}\, [\varphi\, ] = {\L}\, [ z\, (1-M)(1-zM)^{-1}\,{\cal K} \varphi\, ],
\ee
where $M : L_2(m)\to  L_2(m)$ is the multiplication operator 
\be
M\varphi(t) =e^{-t}\varphi(t)
\ee 
and
${\cal K} : L_2(m)\to L_2(m)$ is the integral operator
\be\label{opkappa}
({\cal K} \varphi) (t) =\int_0^\infty  {J_1(2\sqrt{st})\over \sqrt{st}}\,\varphi (s)\,dm(s)\, 
\ee 
and $J_p$ denotes the Bessel function of order $p$.
\end{proposition}
{\sl Proof.}
Letting $f={\L}\, [\varphi\, ]$ we have from (\ref{opevalpose}) and (\ref{doubleid})
\be
{\cal Q}_{z}f (w) =  \sum_{n=1}^\infty {z^n\over (w+n)^2}\int_0^\infty dm(t) \, 
e^{-t/(w+n)}\, \varphi (t)\, .
\ee
Clearly, for $|z|\leq 1$, the sum $\sum_{n=1}^\infty {z^n\over (w+n)^2}e^{-t/(w+n)}$ is uniformly
convergent in $t\in \R^+$. Therefore, interchanging summation and integration we get
\begin{eqnarray}
\sum_{n=1}^\infty {z^n\over (w+n)^2}e^{-t/(w+n)} &=& \sum_{k\geq 0} {(-t)^k\over k!}\sum_{n=1}^\infty
{z^n\over (w+n)^{2+k}}\nonumber \\
&=& \sum_{k\geq 0}{(-t)^k\over k!} \, z \, \Phi (z,k+2,w+1) 
\end{eqnarray}
where $\Phi(z,a,b) = \sum_{n=0}^\infty {z^n\over (b+n)^a}$ is the Lerch transcendental function
which, for $\Re a >1$, possesses the integral representation
\be
\Phi(z,a,b) = \sum_{n=0}^\infty {z^n\over (b+n)^a} = {1\over \Gamma (a)}\int_0^\infty
{s^{a-1}e^{-(b-1)s}\over e^s-z}\, ds\, .
\ee
This yields
\be
\Phi (z,k+2,w+1) = {1\over (k+1)!}\int_0^\infty
{s^{k+1}e^{-ws}\over e^s-z}\, ds\, \cdot
\ee 
Noting that 
\be 
\sum_{k\geq 0} {(-st)^k\over (k+1)!\, k!} = {J_1(2\sqrt{st})\over \sqrt{st}}
\ee
where $J_1(x)$ is the Bessel function of the first kind,
we have thus found that
\begin{eqnarray}\label{repres}
{\cal Q}_{z}f (w) &=& \int_0^\infty ds\,{zs\over (e^s -z)}\,  e^{-ws}\, \int_0^\infty dm(t) \,  {J_1(2\sqrt{st})\over \sqrt{st}}
\varphi (t)\nonumber \\
&=& \int_0^\infty dm(s) \,e^{-ws} \, (z\, (1-M)(1-zM)^{-1}\,{\cal K} \varphi) (s) \\
&=&({\L}\,[z\, (1-M)(1-zM)^{-1}\,{\cal K} \varphi])(w).\nonumber
\end{eqnarray}
Notice that for each $t\in \R^+$ the function 
${J_1(2\sqrt{st})/ \sqrt{st}}$ is uniformly bounded and continuous for $s\in \R^+$. 
It is then an easy task to verify that for $\varphi \in L_2(m)$ and for $|z|\leq 1$ the function $(1-M)(1-zM)^{-1}\,{\cal K} \varphi$
is in $L_2(m)$ as well. $\qed$
\begin{remark}\label{symm}
As already remarked in {\rm \cite{MaR1}}, the integral operator ${\cal K}$ 
is symmetric. Therefore the above Proposition with $z=1$ yields ${\rm sp} \,({\cal Q}) \subset \R$.
\end{remark}
But we can say more. Indeed, the operator $(1-zM)$ is invertible in $L_2(m)$  with bounded inverse
provided $1/z \notin [0,1]$. Therefore,
for any $\varphi \in L_2(m)$ the integral in 
(\ref{rep})
converges uniformly in any compact region of the complex $z$-plane not containing points of the ray
$(1,+\infty)$. Moreover, it has been proved in \cite{MaR1} that the operator ${\cal K}$ is
compact (actually trace-class) in $L_2(m)$. Therefore, as long as $(1-zM)$ has bounded inverse
the operator $(1-M)(1-zM)^{-1}\,{\cal K}$ is compact as well (being the composition of a compact operator with
a bounded operator).
Proposition \ref{invspace} and the above observations prove the following result,
\begin{theorem}\label{analcont}
The operator-valued function $z\to {\cal Q}_{z}$, when acting on 
${\cal H}_1$, can be
analytically continued to the entire $z$-plane with a cross
cut along the ray $(1,+\infty)$, and for each $z$ in this domain is isomorphic 
to the operator 
\be
{\cal K}_{z} :=z\, (1-M)(1-zM)^{-1}\,{\cal K}
\ee 
acting on $L_2(m)$.
They are both compact operators.
\end{theorem}
\begin{remark}\label{noeigenv} The relevance of the above result issues from the following observation: 
the spectral radius of ${\P}$ in any reasonable Banach space of functions
is equal to one (see \cite{Ba1}, \cite{C}) so that according to Proposition \ref{oprel2} 
there are no $z$-values with $|z|<1$ such that 
$1$ is an eigenvalue of ${\cal Q}_z$. Therefore, if we aim to exploit the identities in Proposition \ref{oprelation}
in order to investigate the spectrum of $\P$ 
(when acting upon a suitable function space, see below) it is
necessary to have some analytic continuation of ${\cal Q}_z$ outside the unit disk.
We point out that Proposition \ref{oprelation} and Corollary \ref{oprel2} remain valid when ${\cal Q}_z$
is analytically continued across the cut $(1,+\infty)$.
\end{remark}
\begin{remark}\label{bundi} Putting
\be
H_{\delta} := \{w\in \C : \Re w > \delta\}
\ee
one sees that a function $f={\L}\, [\varphi\, ]$ with $\varphi \in L_2(m)$ can be extended to a function holomorphic in the half-plane $H_{-{1\over 2}}$.

\noindent
If, in addition, $f$ is
an eigenfunction corresponding to a non-zero eigenvalue $\lambda$ of ${\cal Q}_{z}$ in ${\cal H}_1$,
for some non-zero $z\in \C\setminus (1,\infty)$,
then
\be
\lambda\, \varphi (t)= ({\cal K}_{z}\, \varphi)(t) =
 \left({1-e^{-t}\over 1/z-e^{-t}}\right)\,\int_0^\infty  
{J_1(2\sqrt{st})\over \sqrt{st}}\,\varphi (s)\,dm(s)\, .
\ee
Since the integral in the r.h.s. is bounded for all $t\in [0,\infty)$ the function
$\varphi(t)$ is bounded as well is this domain and therefore $f$ 
is holomorphic in the half plane
$H_{-1}$. 
\end{remark}
Putting together the above, Proposition \ref{oprelation} along with standard arguments (see \cite{DS}, Chap. VII) we get,
\begin{corollary}\label{merom}
The operator-valued function $z\to (1-{\cal Q}_{z})^{-1}$, when acting on 
${\cal H}_1$, is analytic in the open unit disk $\{z : |z|<1\}$ and can be
meromorphically continued to the entire $z$-plane with a cross
cut along the ray $[1,+\infty)$. It has a pole
whenever ${\cal K}_{z}$ has $1$ as an eigenvalue. 
\end{corollary}
\vsni
\noindent
Now, from Proposition \ref{oprelation} we obtain the following formal relation for the resolvent
${\cal R}_\lambda$
of ${\P}$:
\be\label{resolvent}
{\cal R}_\lambda \equiv (\lambda-{\P})^{-1}=(\lambda-{{\P}_0})^{-1}(1-{\cal Q}_{1/\lambda})^{-1}.
\ee
The analytic properties of the first factor in the r.h.s. can be 
understood in terms of the spectrum of the operator ${\cal P}_0$ when acting
on a suitable function space invariant under the action of ${\P}$.
A calculation along the same lines as in the
proof of Proposition \ref{invspace} shows that, for $f\in {\cal H}_1$ with $f={\L}[\varphi]$,
\be\label{inj}
(1-z\, {\P}_0)^{-1}f (w)= {1\over w^2} \int_0^\infty e^{-t/w}\,  e^t \,z^{-1}\,({\cal K}_z \varphi)(t)\, dm(t) .
\ee
We shall therefore characterize
the space ${\cal H}_0$ to be acted on by $\P$ as follows:
\begin{definition}
We denote by ${\cal H}_0$ 
the Hilbert space of all complex-valued functions $f$ which can be represented as a generalized 
Borel transform
\be\label{represent2}
f(w) = ({\B}\, [\varphi])(w):={1\over w^2}\int_0^\infty  e^{-t/w}\,e^t\,  \varphi (t)\,dm(t) ,\quad 
\varphi \in L_2(m),
\ee
endowed with the inner product
\be
(f_1,f_2) = \int_0^\infty  {\overline {\varphi_1(t)}}\, \varphi_2  (t)\, dm(t)\quad\hbox{if}\quad f_i={\B}\, [\varphi_i].
\ee
\end{definition}
\begin{remark}\label{another}
A function $f \in {\cal H}_0$ is holomorphic in the disk
\be
D_1=\{w\in \C : \Re {1\over w} > {1\over 2}\}=\{w\in \C : |w-1|<1\}.
\ee
For $w$ real and positive a simple change of variable makes (\ref{represent2}) 
in the form
\be\label{borel}
f(w) = {1\over w}\int_0^\infty  e^{-s}\, \psi (sw)\,ds\quad \hbox{with}\quad
\psi (t) = \left({t\over 1-e^{-t}}\right) \varphi (t)\, .
\ee
\end{remark}
\begin{remark}\label{densities}
The $F$-invariant density $e$ (see (\ref{sigmafinite})) can be
represented as 
\be\label{e}
e=\left({1\over \log 2}\right)\, {\B}\, \left[\, {1-e^{-t}\over t}\,\right],
\ee
whereas for the $G$-invariant density $h$ we have (see also Remark \ref{h})
\be
h=\, \left({1\over \log 2}\right)\, {\L}\, \left[\, {1-e^{-t}\over t}\,\right]=\, \left({1\over \log 2}\right)\,
{\B}\, \left[\, {(1-e^{-t})^2\over t}\,\right].
\ee
In the representation of Remark \ref{another} we have that
if $f=e\cdot \log 2$ then $\psi (t) \equiv 1$ 
whereas for $f=h\cdot \log 2$ we find $\psi (t) =1-e^{-t}$.
Both these functions can be viewed as ordinary Borel transforms of a sequence $\{a_n\}_{n=0}^\infty$, 
i.e. $\psi (t) =\sum_{n=0}^\infty t^n a_n /n!$ so that by (\ref{borel}) we have 
$w\cdot f(w)= \sum_{n=0}^\infty w^n a_n$.
In the former case we find $a_0=1$ and $a_n =0$ for $n>0$,
in the latter $a_0=0$ and $a_n=(-1)^{n-1}$ for $n>0$. Therefore in both cases the integral (\ref{borel}) provides a continuation of 
$w\cdot f(w)$ outside the disk $D_1$ (see \cite{Tit1}, p.164).
\end{remark}
We now have the following,
\begin{lemma} \label{lemmino}
For all $\varphi \in L_2(m)$
\be 
{\cal L}\, [\,\varphi] = {\cal B}\, [\,(1-M)\,{\cal K}\, \varphi]
\ee
where $M\varphi(t) =e^{-t}\varphi(t)$ and ${\cal K}$ is the symmetric integral operator  
defined in (\ref{opkappa}).
\end{lemma}
{\sl Proof.} The proof is an easy calculation based on Tricomi's theorem (see \cite{Sne}, p.165)
\be
{1\over u^{p+1}}\int_0^\infty dt \, e^{-t/u} \varphi (t) = \int_0^\infty dt \, e^{-tu}\,
\int_0^\infty ds\, \left({t\over s}\right)^{p\over 2}\, J_p(2\sqrt{st})\, \varphi (s),
\ee
with $p=1$, and therefore we omit it. $\qed$
\vsni
\noindent
It is now not difficult to verify that
\be\label{piuno}
{\P}_1\, {\B} [\varphi]  = {\L}\, [\varphi],
\ee
and
\be\label{pizero}
{\P}_0\, {\B} [\varphi]= {\B}\, [\, M\varphi\,].
\ee
In addition we have
\be\label{shift}
{\cal S}{\L}[\varphi] = {\L}\, [M\varphi],
\ee
so that
\be
{\P}_1\,{\P}_0^{n-1}\, {\B} [\varphi]={\cal S}^{n-1}\,  {\P}_1 {\B} [\varphi]={\L}\, [\,M^{n-1} \varphi\,],
\ee
and therefore
\be \label{emmezeta}
{\cal Q}_z {\B}\, [\,\varphi\,] = z\cdot {\L}\, [\,(1-zM)^{-1} \varphi\,].
\ee
We thus see that ${\P}_0$ leaves ${\cal H}_0$ invariant and by (\ref{shift}) its spectral properties in
${\cal H}_0$ are identical to those of ${\cal S}$ in ${\cal H}_1$. Moreover ${\P}_1$ maps ${\cal H}_0$ into
${\cal H}_1 \subset {\cal H}_0$, and the same does ${\cal Q}_z$ for all $z\in \C \setminus (1,+\infty)$.
Notice that using Lemma \ref{lemmino} and (\ref{emmezeta}) we immediately recover Proposition \ref{invspace},
in that
\be
{\cal Q}_z {\L}\, [\, \varphi\,] = {\cal Q}_z {\cal B}\, [\,(1-M)\,{\cal K}\, \varphi] = 
 {\L}\, [\,z\cdot (1-zM)^{-1} (1-M)\,{\cal K}\, \varphi\,] \equiv {\L}\, [ \,{\cal K}_z \varphi\, ].
\ee
\noindent
We are now in the position to write explicit representations for $\P$ and its resolvent 
${\cal R}_\lambda$ in the space
${\cal H}_0$.
\begin{theorem} \label{resrep}
Let $f\in {\cal H}_0$, that is $f={\B}\, [\varphi]$ for some $\varphi \in L_2(m)$, then
\be\label{pi}
{\P} f = {\B}\, [\, (M+ (1-M){\cal K}\,)\varphi\,],
\ee
and
\be\label{resol}
{\cal R}_\lambda f \equiv (\lambda -{\P})^{-1} f\, = 
{\B}\, [\, (1-{\cal K}_{1/\lambda})^{-1}(\lambda-M)^{-1}\varphi \, ].
\ee
\end{theorem}
\begin{remark}\label{analyticityimproving}
Note that for $\varphi \in L_2(m)$ the functions 
\be 
M\,\varphi\quad\hbox{and}\quad (1-M){\cal K}\,\varphi
\ee
are bounded at infinity and therefore, by (\ref{pi}), the function ${\P} f$ with
$f={\B}\, [\varphi]$ is analytic in the half-plane $H_0$. In particular so is any eigenfunction of
$\P$ in ${\cal H}_0$.
\end{remark}
{\sl Proof of Theorem \ref{resrep}}. From
(\ref{pizero}) and (\ref{piuno}) one obtains ${\P} f =  {\B}\, [\, M\varphi\,] + {\L}\, [\varphi]$,
so that (\ref{pi}) follows using Lemma \ref{lemmino}.
The expression for ${\cal R}_\lambda$ can now be obtained directly from (\ref{pi}). But we can also
make use of (\ref{emmezeta}) and  (\ref{repres}) to obtain, for a given $f={\B}\, [\varphi]$,
\be
{\cal Q}^n_{1/\lambda} \, f = {\L}\, [\, {\cal K}_{1/\lambda}^{n-1} (\lambda-M)^{-1}\varphi\,]
\ee
and therefore
\be 
(1-{\cal Q}_{1/\lambda})^{-1}f = {\B}\, [\varphi] + {\L}\, [\,(1- {\cal K}_{1/\lambda})^{-1}\, (\lambda-M)^{-1}\varphi\,].
\ee
This expression along with (\ref{resolvent}), (\ref{inj}) and (\ref{pizero})
yield
\begin{eqnarray}
{\cal R}_\lambda f &=& {\B}\,[\, (\lambda-M)^{-1}\, \varphi\,]+
{\B}\,[\,{\cal K}_{1/\lambda}(1- {\cal K}_{1/\lambda})^{-1}\, (\lambda-M)^{-1}\,\varphi\,]\nonumber \\
&=& {\B}\, [\, (1-{\cal K}_{1/\lambda})^{-1}(\lambda-M)^{-1}\varphi \, ].\nonumber 
\end{eqnarray}
Using Corollary \ref{merom} we see that ${\cal R}_\lambda$ extends to a meromorphic (operator-valued) function in ${\overline \C} \setminus [0,1]$. $\qed$
\vsni
\noindent
The next theorem (partially) describes the spectrum of $\P$ in ${\cal H}_0$.
\begin{theorem} \label{spectrump}
The spectrum of the operator
${\P}:{\cal H}_0 \to {\cal H}_0$ is the union of
$[0,1]$ and a finite or countably infinite set of eigenvalues of finite multiplicity.
\end{theorem}
{\sl Proof.}
By Theorem
\ref{resrep} the action of transfer operator $\P$ on ${\cal H}_0$ can be explicitly expressed in the form
\be
{\P}{\B}\, [\varphi] = {\B}\, [\, T\varphi\,],
\ee
with 
\be\label{sum}
(T\varphi)(t):= e^{-t}\varphi (t)+ 
\int_0^\infty K(s,t)\varphi (s) ds 
\ee
and
\be  
K(s,t) = e^{-t}\, \left( {e^t-1\over e^s-1}\right)\, \sqrt{s\over t}\,\, J_1(2\sqrt{st}).
\ee
It
is an easy exercise to check that $M$ when acting upon $L_2(m)$ is self-adjoint and its spectrum is the line segment 
$[0,1]={\rm Cl}\, \{ e^{-t}\, :\, t\in \R^+\}$ (see, e.g., \cite{DeV}). 
Therefore the spectrum of $\P$ in ${\cal H}_0$ is given by a compact perturbation of the continuous
spectrum $\sigma_c = [0,1]$.
The assertion is now a consequence of Theorem 5.2 in \cite{GK}.
$\qed$

\vsni
\noindent
We shall now characterize some properties of the eigenfuctions of $\P$ in ${\cal H}_0$. 
First, it is  easy to see that $\lambda =0$ is an eigenvalue of infinite
multiplicity. This follows by noting that
(see (\ref{invbran}) and (\ref{opP})) any 
function $f\in {\cal H}_0$ which is odd w.r.t. $x=1/2$, e.g. 
$f(w)=1-2w = {\B}\, [\, (1-t)(1-e^{-t})\,]$ lies in the kernel of $\P$.

\noindent
Now suppose that ${\P}f = \lambda \,f$ for some $f\in {\cal H}_0$ and $\lambda \ne 0$, or
explicitly
\be\label{eigeneq}
\lambda \, f (w) = \left({1\over w+1}\right)^2 \left[ f\left( {w\over w+1}\right) + 
f\left( {1\over w+1}\right) \right].
\ee
By Remark \ref{analyticityimproving} $f(w)$ extends analytically to the half-plane $H_0$.
If we transform this equation by substituting $1/w$ for $w$ and then dividing through $w^2$ we get
\be\label{eigeneq1}
\lambda \, w^{-2}\, f\left({1\over w}\right) = \left({1\over w+1}\right)^2 \left[ f\left( {1\over w+1}\right) + 
f\left( {w\over w+1}\right) \right].
\ee
Therefore $f$ satisfies
\be\label{funceq}
w\, f(w) = {1\over w}\, f\left({1\over w}\right)
\ee
for all $w\in H_0$.
Note that applying (\ref{funceq}) to each term of the r.h.s. in (\ref{eigeneq}) one obtains
\be\label{dudu}
\lambda \, w\, f(w) = w\, f(w+1) + {1\over w}\, f\left(1+{1\over w}\right).
\ee
For $\lambda =1$ this yields $w\, f(w)=1$.
Note that for $f ={\B}\, [\varphi]$ we have
\be\label{lolo}
w^{-2}\, f\left({1\over w}\right) =\int_0^\infty  e^{-t\, w}\,e^t\,  \varphi (t)\,dm(t) = 
{\B}\, [(1-M){\cal K}\,M^{-1}\varphi].
\ee
Therefore the functional equation (\ref{funceq}) can be written as
\be\label{ago}
(1-M){\cal K}\,M^{-1}\varphi = \varphi. 
\ee
Now, given a continuous function
$\psi$ on $\R^+$ one can define (a version of) its {\sl Hankel transform} (of order $1$) as the integral
\be\label{hankel}
({\cal J} \psi) (t) = \int_0^\infty J_1(2\sqrt{st})\,\sqrt{t\over s}\, \psi (s)\, ds.
\ee
From the estimates $J_1(t) \sim t$ as $t\to 0^+$ and $J_1(t)= O(t^{-1/2})$ 
as $t\to \infty$ (\cite{E}, vol.II) we see that
the conditions on $\psi$ sufficient to give the absolute convergence of the integral (\ref{hankel}) are
$\psi (t) = O(t^ {-\beta})$ as $t\to \infty$ with $\beta > -1/4$ and 
$\psi (t) = O(t^ {\alpha})$ as $t\to 0^ +$ with $\alpha > -1$.
The identity (\ref{ago}) then says that the function (cf. Remark \ref{another})
\be\label{lilla}
\psi (t) = \left({t\over 1-e^{-t}}\right) \varphi (t)
\ee  
 satisfies 
\be\label{equa}
 \psi (t) = \int_0^\infty J_1(2\sqrt{st})\,\sqrt{t\over s}\, \psi (s)\, ds.
\ee
Note that the simplest solution of this equation is $\psi \equiv 1$ and corresponds to $f=e$
(more general self-reciprocal functions satisfying equations related to (\ref{equa}) are discussed, e.g., in the book \cite{Tit2}).
Furthermore, putting together (\ref{funceq}), (\ref{lolo}) and (\ref{lilla}) we have that 
\be\label{olt}
f(w) = \int_0^\infty  e^{-t\, w}  \psi (t)\,dt
\ee
for all $w\in H_0$. Finally, one easily checks that if $\varphi \in  L_2(m)$ then $\psi \in  L_2({\hat m})$ where 
$$d{\hat m}(t)={e^{-t}(1-e^{-t})\over t\, \cdot \, \log{2}}\, dt$$ We summarize the above in the following
\begin{theorem} \label{hank}
If $f\in {\cal H}_0$ satisfies  ${\P}f = \lambda \,f$ for some $\lambda \ne 0$
then $f$ is the (ordinary) Laplace transform of a function $\psi \in  L_2({\hat m})$ 
which is self-reciprocal w.r.t. Hankel transform of order $1$, 
namely $f$ and $\psi$ satisfy (\ref{olt}) and (\ref{equa}), respectively.
\end{theorem}

\noindent
Now from Corollary \ref{oprel2} we know that a function 
$f={\B}\, [\varphi]$ satisfies ${\P}f = \lambda \,f$
if and only if (the analytic continuation of)
${\cal K}_{1/\lambda}:L_2(m)\to L_2(m)$ satisfies ${\cal K}_{1/\lambda}\,\varphi =\varphi$, which
can also be written as
\be\label{dudu}
({\cal K} \, \varphi) (t) = {\lambda -e^{-t}\over 1-e^{-t}}\, \varphi(t)= {\lambda -e^{-t}\over t}\, \psi(t)\ .
\ee
Expressing the integral operator ${\cal K}$ in terms of the Hankel transform (\ref{hankel}) we get
$({\cal K} \, \varphi)(t) ={1\over t}\, {\cal J} (\exp_{-1}\cdot \psi)(t)$, where we have defined the 
function $\exp_c : \R \to \R$ by $\exp_c (t) =e^{ct}$. Identities  (\ref{equa}) and (\ref{dudu}) then yield the integral equation
\be\label{equona}
{\cal J} (\exp_{-1}\cdot \psi) = (\lambda -\exp_{-1})\cdot {\cal J} \psi\, .
\ee
Once more, $\psi \equiv 1$ satisfies this equation with $\lambda =1$ 
(recall that ${\cal J} \exp_{-1}=1-\exp_{-1}$). On the other hand, the above discussion suggests that there are no $\lambda \in \C\setminus \{0,1\}$ such that 
(\ref{equona}) has a (non-constant) solution $\psi \in  L_2({\hat m})$. We are thus are led to formulate the following,
\begin{conjecture} 
The only (non-zero) eigenvalue of 
${\P}:{\cal H}_0 \to {\cal H}_0$ is $\lambda =1$.
\end{conjecture}

\vsni
\noindent 
We end this Section with two additional remarks.
\begin{remark} 
(\ref{dudu}) is a particular case of the {\rm Lewis functional equation} 
\be
f (w) -f(w+1)={1\over w^{2(q+1)}} f\left( 1+{1\over w}\right),
\ee
which is related to the so called Maass cusp forms, i.e. PSL($2,\Z$)-invariant eigenfunctions 
of the Laplacian on the Poincar\'e upper half-plane
which vanish at the cusp (see \cite{Le}). Another type of functions equivalent to (even) 
Maass forms and considered in \cite{Le} are those satisfying an integral equation which in our notation
writes
\be
g(t) = \int_0^\infty \,  {J_{2q+1}(2\sqrt{st})\over \sqrt{st}} \left({s\over t}\right)^{q}
 \, g (s) \, dm(s)\,.
\ee
By the foregoing (see Remark \ref{densities}) we see that for $q=0$ we have the relation
\be
f= {\cal B}\, [g].
\ee
\end{remark}
\begin{remark}
In the recent work \cite{P2}, following \cite{P1} ten years later and somehow inspired by the construction presented here, Thomas Prellberg 
has studied the spectrum of (a generalized version of) ${\cal P}$ in a space of functions which is identical
to ${\cal H}_0$ with the exception that the measure on $\R^+$ is slightly different from (\ref{measure}), being given by
\be\label{newmeasure}
d{\tilde m}(t) =t\, e^{-t} \, dt.
\ee
It is easy to see that 
with this new measure the operator ${\cal Q}_z$ is isomorphic under generalized Laplace transform (cf. Theorem \ref{analcont}) to
${\tilde {\cal K}}_z : L_2({\tilde m})\to L_2({\tilde m})$ given by 
\be
{\tilde {\cal K}}_z = z\, (1-zM)^{-1}{\tilde {\cal K}},
\ee
where 
\be
{\tilde {\cal K}} \varphi (s) =  \int_0^\infty d{\tilde m}(t)\, {J_1(2\sqrt{st})\over \sqrt{st}}\, \cdot
\varphi (t) 
\ee
Notice that ${\cal K}_1= (1-M)^{-1}{\tilde {\cal K}}$ which is not symmetric anymore (cf. Remark \ref{symm}).
On the other hand,
the relation given by Lemma \ref{lemmino} now writes (we keep using the symbols ${\cal L}$ and ${\cal B}$
to denote generalized Laplace and Borel transforms w.r.t. the measure ${\tilde m}$):
\be 
{\cal L}\, [\,\varphi] = {\cal B}\, [\,{\tilde {\cal K}}\, \varphi]
\ee
and hence the integral representation of ${\cal P}$ becomes
\be
{\P} {\cal B}\, [\,\, \varphi]= {\B}\, [\, (M+ {\tilde {\cal K}}\,)\varphi\,],
\ee
which is now symmetric (cf. (\ref{pi})). Thus, everything goes as if the operators ${\cal P}$ and ${\cal Q}$
were not `symmetrizable' both at the same time. Also notice that the function $\log 2 \cdot e$
if expressed as a generalized Borel transform now yields the function $\varphi (s) =1/s$
which is not in $L_2({\tilde m})$. 
\end{remark}

\section{Zeta functions}\label{ZF}

We now consider the dynamical zeta functions $\zeta_F$ and $\zeta_G$
associated to the maps $F$ and $G$, respectively, and defined by 
the following formal series \cite{Ba2}:
\be
\zeta_F (z) = \exp \sum_{n=1}^{\infty} {z^n\over n} Z_n(F) \quad\hbox{and}\quad
\zeta_G (s) = \exp \sum_{n=1}^{\infty} {s^n\over n} Z_n(G),
\ee
where the `partition functions' $Z_n(F)$ and $Z_n(G)$ are given by
\be
Z_n(F) =\sum_{x=F^n(x)} \prod_{k=0}^{n-1}{1\over |F'(F^k(x))|}
\quad\hbox{and}\quad
Z_n(G) =\sum_{x=G^n(x)} \prod_{k=0}^{n-1}{1\over |G'(G^k(x))|}\cdot
\ee
Let us first examine how $\zeta_F (z)$ and $\zeta_G (z)$ are related
to one another. Let 
${\rm Per}\, F$ and ${\rm Per}\, G$ denote the sets of all periodic points
of the maps $F$ and $G$, respectively. It is not difficult to realize that,
as subsets of $\ui$, ${\rm Per} \,F\setminus \{0\}={\rm Per}\, G$.
Accordingly, given $x$ in either of these sets, we let $p_F(x)$ and $p_G(x)$
denote the periods of $x$ w.r.t. to $F$ and $G$, respectively. 
They are related by
\be
p_F(x) = \tau (x) + \tau (G(x)) + \cdots + 
\tau(G^{\, p_G(x)-1}(x))\, \cdot 
\ee
Moreover from the definitions of $F$ and $G$ we have
\be
\prod_{k=0}^{p_F(x)-1}{1\over |F'(F^k(x))|}=
\prod_{k=0}^{p_G(x)-1}{1\over |G'(G^k(x))|}= \prod_{k=0}^{p_G(x)-1}\, (G^k(x))^2\, \cdot
\ee
Using this facts we write $Z_n(F)$ as follows:
\be
Z_n(F) = 1+ \sum_{m=1}^n {n\over m}
\sum_{\scriptstyle x=F^n(x)=G^m(x)}
\prod_{k=0}^{m-1}\, (G^k(x))^2\, \cdot
\ee
The second sum ranges over the $n-1 \choose m-1$ ways to write the integer $n$
as a sum of $m$ positive integers.
Therefore, 
\begin{eqnarray}
\sum_{n=1}^{\infty}{z^n\over n} Z_n(F) &=& \log \left({1\over 1- z}\right) +
\sum_{n=1}^{\infty}\sum_{m=1}^{n}
{1\over m} \sum_{\scriptstyle x=F^n(x)=G^m(x)} 
z^n \prod_{k=0}^{m-1}\, (G^k(x))^2\, \nonumber \\
&=& \log \left({1\over 1- z}\right) +
\sum_{\ell=1}^{\infty} {1\over \ell} \sum_{\scriptstyle x=G^\ell(x) }
 z^{p_F(x)}\prod_{k=0}^{\ell-1}\, (G^k(x))^2\,\cdot \nonumber
\end{eqnarray}
We are thus led to study the `grand partition function' $\Xi_\ell(z)$ given by
\be\label{Xi}
\Xi_\ell(z) := \sum_{\scriptstyle x=G^\ell(x) }
z^{p_F(x)}\prod_{k=0}^{\ell-1}\, (G^k(x))^2\,
= \sum_{n=0}^{\infty}z^{\ell+n}\sum_{\scriptstyle x=G^\ell(x)=F^{\ell+n}(x)}
 \prod_{k=0}^{\ell-1}(G^k(x))^2.
\ee
The sum over periodic points yields  
${n+\ell -1 \choose \ell -1}={n+\ell -1 \choose n}$ terms, corresponding to the number of
ways of distributing $n$ identical objects into $\ell$ distinct boxes.
According to (\ref{Gaction}), (\ref{Faction}) and (\ref{Xi}) we can also write $\Xi_\ell(z)$ in the following way:
\be
\Xi_\ell(z)  = \sum_{n=0}^{\infty}z^{\ell+n}\sum_{k_1+\dots +k_\ell=n+\ell}\,  
\prod_{i=1}^\ell x^{2}_{k_i\dots k_\ell k_1\dots k_{i -1}} ,
\ee
where $x_{k_1\dots k_\ell}=[{\overline {k_1,\dots, k_\ell}}]$ denotes the irrational number whose continued fraction expansion is periodic
of period $\ell$ and starts with the entries $k_1,\dots ,k_\ell$.
\noindent
Putting together the above observations we obtain the next result, to be compared with Proposition \ref{oprelation}:
\begin{proposition}\label{twozeta} Consider the 
two-variable zeta function given by
\be\label{twovar}
\z_2 (s,z) = \exp \sum_{\ell=1}^{\infty} {s^\ell\over \ell}\, \Xi_\ell(z).
\ee
Then
we have:
\be
\z_2 (1,z) = (1- z)\, \z_F (z) \quad\hbox{and}\quad \z_2 (s,1) = \z_G (s)
\ee
wherever the series expansions converge absolutely.
\end{proposition}
In order to study the analytic properties of the function $\z_2 (s,z)$ 
we further generalize (\ref{opevalpose}) by introducing a family of operator-valued functions ${\cal Q}_{z,q}$,
$q=0,1,\dots $,
acting as
(see \cite{Ma1} and \cite{D} for related quantities)
\be
{\cal Q}_{z,q} f(x) =(-1)^q\sum_{n=1}^{\infty} z^{\tau(\Phi_n (x))} \cdot|\Phi'_n(x)|^{1+q}\cdot f (\Phi_n (x)),
\ee
together with a family of function spaces ${\cal H}_{1,q}\subseteq {\cal H}_{1}$ such that 
a function $f\in {\cal H}_{1,q}$ can be represented as
\be\label{represent1}
f(w) =({\cal L}_q\, [\varphi])(w):=\int_0^\infty dm(t) \, e^{-tw}\,t^{q}\, \varphi (t), \quad 
\varphi \in L_2(m) .
\ee
In particular ${\cal Q}_{z,0}\equiv {\cal Q}_z$, ${\cal L}_0\equiv {\cal L}$ and ${\cal H}_{1,0}\equiv {\cal H}_{1}$.
We have the following result.
\begin{proposition} For any given $q=0,1\dots$ the operator valued function $z\to {\cal Q}_{z,q}$
when acting on ${\cal H}_{1,q}$
can be
analytically continued to the entire $z$-plane with a cross
cut along the ray $(1,+\infty)$. For each $z$ in this domain we have
\be
{\cal Q}_{z,q}\,{\L}_q\, [\varphi\, ] = {\L}_q\, [ \,{\cal K}_{z,q} \varphi\, ],
\ee
where ${\cal K}_{z,q} : L_2(m)\to L_2(m)$ is given by
\be\label{opj}
({\cal K}_{z,q} \varphi) (t) :=  (-1)^q\, z\, (1-M)(1-zM)^{-1} \,\int_0^\infty dm(s)\, 
{J_{2q+1}(2\sqrt{st})\over \sqrt{st}}
 \, \varphi (s) \, .
\ee 
The operators ${\cal Q}_{z,q}: {\cal H}_{1,q} \to {\cal H}_{1,q}$ and ${\cal K}_{z,q} : L_2(m)\to L_2(m)$ are 
both of trace class. 
\end{proposition}
{\sl Proof.} The first part follows from a straightforward extension to non zero $q$ values of the arguments of the
previous Section. The proof of the last assertion can be extracted from (\cite{Ma1}, Theorem 3). $\qed$

\vsni
\noindent
Now, the trace of the operator ${\cal K}_{z,q}$ is easily obtained (see also \cite{Ma1}):
\begin{eqnarray}\label{trace}
{\rm tr}\,\, {\cal K}_{z,q} &=& (-1)^q\,z\, \int_0^\infty \,
{J_{2q+1}(2t)\over e^t-z}\, dt \nonumber \\
&=& (-1)^q\,\sum_{k=1}^\infty\, z^k \, \int_0^\infty e^{-kt}\,J_{2q+1}(2t)\, dt  \\
&=& (-1)^q\,\sum_{k=1}^\infty\, z^k\, {\;\; x_k^{2(q+1)}\over 1+x_k^2},
\nonumber
\end{eqnarray}
where the numbers $x_k ={\sqrt{k^2+4}-k\over 2}=[k,k,k,\dots ]\equiv [{\overline k}]$ are the fixed points of $G(x)$ and
the identity \cite{GR}
\be
\int_0^\infty e^{-kt}J_p(2t)dt = {(\sqrt{k^2+4}-k)^p\over 2^p\, \sqrt{k^2+4}},\quad p=0,1,\dots 
\ee
has been used. From (\ref{trace}) we immediately obtain the trace formula
\be
\Xi_1(z) = {\rm tr}\,\, {\cal K}_{z,0} - {\rm tr}\,\, {\cal K}_{z,1}.
\ee
But we can say more. Indeed, a straightforward adaptation of (\cite{Ma1}, Corollaries 4 and 5) to our 
$z$-dependent situation leads to 
the following general expressions:
\be\label{traceformula}
\Xi_\ell(z) = {\rm tr}\,\, {\cal K}_{z,0}^\ell - {\rm tr}\,\, {\cal K}_{z,1}^\ell
={\rm tr}\,\, {\cal M}_{z,0}^\ell - {\rm tr}\,\, {\cal M}_{z,1}^\ell,
\ee
with
\be
{\rm tr}\,\, {\cal K}^\ell_{z,q}=(-1)^{q\ell}\,\sum_{k_1,\dots ,k_\ell=1}^\infty\, z^{k_1+\cdots +k_\ell}
\,{ \prod_{i=1}^\ell x^{2(q+1)}_{k_i\dots k_\ell k_1\dots k_{\ell -1}}\;\; \over 1-(-1)^\ell 
\prod_{i=1}^\ell x^{2}_{k_i\dots k_\ell k_1\dots k_{\ell -1}} }\, \cdot
\ee
Formula (\ref{traceformula}) along with standard arguments (see \cite{Ma1}) allow us to write the two-variables
zeta function (\ref{twovar}) as a ratio of Fredholm determinants,
\be
\z_2 (s,z) = \exp \sum_{\ell=1}^{\infty} {s^\ell\over \ell}\, \Xi_\ell(z)
={{\rm det}\, (1-s\,{\cal K}_{z,1}) 
\over {\rm det}\, (1-s\,{\cal K}_{z,0}) } = {{\rm det}\, (1-s\,{\cal M}_{z,1}) 
\over {\rm det}\, (1-s\,{\cal M}_{z,0}) }\, ,
\ee
where by definition 
\be
{\rm det}\,(1-s\,{\cal K}_{z,q}) = \exp \left( -\sum_{\ell =1}^\infty {s^\ell\over \ell}\, {\rm tr}\,\, 
{\cal K}^\ell_{z,q}
\right)
\ee
is in the sense of Grothendieck \cite{G}.
We have thus proved the following result. 

\begin{theorem}\label{twoanal}
Set ${\cal K}_{z}\equiv {\cal K}_{z,0}$, then we have:
\noindent
\begin{enumerate}
\item for each $s\in \C$, the function $\z_2 (s,z)$, considered as a
function of the variable $z$,
extends to a meromorphic function in the
cut plane 
$\C\setminus [1,\infty)$. Its poles are located among 
those $z$-values such
that ${\cal K}_{z}:L_2(m)\to L_2(m)$  has $1/s$ as an eigenvalue;
\item for each $z\in \C \setminus (1,\infty)$, the function $\z_2 (s,z)$, considered as a
function of the variable $s$,
extends to a meromorphic function in  $\C$. Its poles are located among the inverses of the
eigenvalues of ${\cal K}_{z}:L_2(m)\to L_2(m)$.
\end{enumerate}
\end{theorem}
Putting together the above Theorem and Proposition \ref{twozeta} we obtain 
\begin{corollary} \label{zetas}
The dynamical zeta functions $\zeta_F$ and $\zeta_G$ of the Farey and Gauss maps
have the following properties:
\begin{enumerate}
\item $\zeta_F (z)$ has a meromorphic extension to the cut plane $\C \setminus [1,\infty)$;
\item  $\zeta_G (s)$ has a meromorphic extension to $\C$.
All poles are real and are located among the inverses of
the eigenvalues of ${\cal K}:L_2(m)\to L_2(m)$. 
\end{enumerate}
\end{corollary}
\begin{remark} Statement 1 of Corollary \ref{zetas} is akin to Corollary 1.3 in \cite{Rug}. On the other hand, the validity of Conjecture 1 would imply that
 $\zeta_F (z)$ is actually analytic in $\C \setminus [1,\infty)$.
Statement 2 was proved by Mayer in
\cite{Ma3}, where he also showed (using results from \cite{Rue}) 
that the poles of $\zeta_G (s)$, if arranged in increasing
absolute values and according to their order tend to infinity exponentially fast.
\end{remark}


\begin{thebibliography}{DEGHL}


\bibitem[Ba1]{Ba1}
{\sc V Baladi}, \, {\sl Positive Transfer Operators and Decay of Correlations}, World Scientific, 2000. 

\bibitem[Ba2]{Ba2}
{\sc V Baladi}, \,{\sl Dynamical zeta functions}, Real and Complex
Dynamical Systems (B. Branner and P. Hjorth eds.), Kluwer Academic
Publishers, 1997. 

\bibitem[Bi]{Bi}
{\sc P Billingsley}, \, {\sl Ergodic Theory and Information}, John Wiley, 1964. 

\bibitem[C]{C}
{\sc P Collet}, {\sl Some ergodic properties of maps of the interval},  
in Dynamical Systems, proceedings of the first UNESCO CIMPA school of Dynamics and disordered systems
(Temuco, Chile, 1991), Herman (1996).

\bibitem[CI]{CI} {\sc P Collet, S Isola},
{\it  On the Essential Spectrum of the Transfer Operator for Expanding Markov Maps},  Commun. Math. Phys. {\bf 139} 
(1991), 551-557.


\bibitem[DeV]{DeV}
{\sc C L DeVito}, {\sl Functional Analysis and Linear Operator Theory},  Addison-Wesley Publ. Co. 1990.

\bibitem[D]{D} {\sc P Dodds}, Master's thesis, The University of Melbourne  (1993).

\bibitem[DS]{DS}
{\sc N Dunford, J T Schwartz}, {\sl Linear Operators, Part I}, Interscience Publ., Inc., New York, 1963.


\bibitem[GR]{GR}
{\sc I Gradshteyn, I Ryzhik}, {\sl Table of integrals, series and products}, Academic Press, 1965.

\bibitem[E]{E}
{\sc A Erd\`ely et al.}, {\sl Higher trascendental functions} (Bateman manuscript project), 
Vols. I-III, McGraw-Hill, New York, 1953-1955.

\bibitem[F]{F} {\sc M Feigenbaum},
{\it Parabolic rational maps}, J. Stat. Phys. {\bf 52} (3-4) (1988), 527-569.


\bibitem[G]{G}
{\sc A Grothendieck}, {\it La theorie de Fredholm}, Bull. Soc. Math. Fr., 
{\bf 84} (1956), 319-384.

\bibitem[GK]{GK}
{\sc G C Goheberg, M G Krejn}, {\sl Introduction \`a la th\'eorie des op\'erateurs lin\'eaires non auto-adjoints
dans un espace Hilbertien}, Dunod, Paris, 1971.


\bibitem[HR]{HR} {\sc G H Hardy, E M Wright},
{\sl An Introduction to the Theory of Numbers},  Clarendon Press, Oxford, 1979.


\bibitem[HI]{HI} {\sc N T A Haydn, S Isola},
{\it Parabolic rational maps}, J. London Math. Soc. {\bf 63} (2) (2001), 673-689.

\bibitem[HY]{HY}  {\sc H Hu, L-S Young}, {\it Nonexistence of SRB measures for 
some systems that are ``almost Anosov"},  Erg. Th. Dyn. Syst. {\bf
15} (1995), 67-75.


\bibitem[Is]{Is} {\sc S Isola},
{\it On systems with finite ergodic degree}, 
Preprint.

\bibitem[Ki]{Ki} {\sc A Khinchin},
{\it Continued Fractions}, 
University of Chicago Press, 1964.


\bibitem[Le]{Le}
{\sc J B Lewis}, {\it Spaces of holomorphic functions equivalent to the even Maass cusp forms}, 
Invent. Math. {\bf 127} (1997), 271-306. 

\bibitem[Ma1]{Ma1} {\sc D H Mayer}, {\it On the thermodynamic formalism for the Gauss map}, 
Commun. Math. Phys. {\bf 130} (1990), 311-333.


\bibitem[Ma2]{Ma2} {\sc D H Mayer}, {\it Continued fractions and related transformations}, in
Ergodic Theory, Symbolic Dynamics and  Hyperbolic Spaces, T Bedford, M Keane and C Series Eds., 
Oxford University Press, 1991.

\bibitem[Ma3]{Ma3} {\sc D H Mayer}, {\it On a $\zeta$ function related to the continued fraction transformation}, 
Bull. Soc. math. France {\bf 104} (1976), 195-203.


\bibitem[MaR1]{MaR1} {\sc D H Mayer, G Roepstorff}, {\it On the relaxation time of Gauss' continued fraction map I.
The Hilbert space approach.}, 
J. Stat. Phys. {\bf 47} (1987),149-171.

\bibitem[MaR2]{MaR2} {\sc D H Mayer, G Roepstorff}, 
{\it On the relaxation time of Gauss' continued fraction map II.
The Banach space approach.}, 
J. Stat. Phys. {\bf 50} (1988),331-344.

\bibitem[Me]{Me} {\sc C Meunier }, {\it Continuity of type I intermittency from a
measure theoretical point of view}, J. Stat. Phys. {\bf 36} (1984), 321-365.

\bibitem[PM]{PM} {\sc Y Pomeau, P Manneville }, {\it Intermittent transition to
turbulence in dissipative dynamical systems}, 
Comm. Math. Phys. {\bf 74} (1980), 189-197.

\bibitem[P1]{P1} {\sc T Prellberg}, Ph.D. thesis, Virginia Tech (1991).

\bibitem[P2]{P2} {\sc T Prellberg}, Preprint (2001).

\bibitem[PS]{PS} {\sc T Prellberg, J Slawny}, {\sl Maps of intervals with 
indifferent fixed points: thermodynamic formalism and phase transitions},
J. Stat. Phys. {\bf 66} (1992), 503-514.


\bibitem[Rue]{Rue}
{\sc D Ruelle}, \,{\it Zeta functions for expanding maps and Anosov flows}, 
 Invent. Math. {\bf 34} (1976), 231-242.

\bibitem[Rug]{Rug}
{\sc H H Rugh}, {\it Intermittency and Regularized Fredholm Determinants}, 
Invent. Math. {\bf 135} (1999), 1-24.

\bibitem[Sne]{Sne}
{\sc I N Sneddon}, {\it The use of Integral Transforms}, 
Tata McGraw-Hill Publ. Co. Ltd, New Delhi, 1974.


\bibitem[Tit1]{Tit1}
{\sc E C Titchmarsh}, {\it The Theory of Functions}, 
Oxford University Press, 1939.

\bibitem[Tit2]{Tit2}
{\sc E C Titchmarsh}, {\it Introduction to the Theory of Fourier Integrals}, 
Oxford at the Clarendon Press, 1937.

\end{thebibliography}
\end{document}